\documentclass{amsart}
\usepackage{latexsym}
\usepackage{amsfonts}
\usepackage{amssymb}
\usepackage{graphicx}
\usepackage{graphicx}

\newtheorem{theorem}{Theorem}[section]
\newtheorem{corollary}[theorem]{Corollary}
\newtheorem{lemma}[theorem]{Lemma}

\theoremstyle{definition}

\theoremstyle{remark}
\newtheorem{remark}{Remark}[section]
\newtheorem*{example}{Example}

\newcommand{\g}{\mathfrak{g}}

\newcommand{\ak}{\mathfrak{k}}
\newcommand{\ap}{\mathfrak{p}}

\newcommand{\RR}{\mathbb{R}}

\newcommand{\diag}{\mathrm{diag}}

\newcommand{\Ad}{\mathrm{Ad}}
\newcommand{\ad}{\mathrm{ad}}
\newcommand{\im}{\mathrm{Im}}

\begin{document}

\title[]
{On the realization of Riemannian symmetric spaces in Lie
groups$^1$}

\author[]{Jinpeng An}
\address{School of mathematical science, Peking University, Beijing, 100871, P. R. China }
\email{anjinpeng@pku.edu.cn}

\author[]{Zhengdong Wang}
\address{School of mathematical science, Peking University, Beijing, 100871, P. R. China}
\email{zdwang@pku.edu.cn}

\maketitle

{\bf Abstract.} In this paper we give a realization of some
symmetric space $G/K$ as a closed submanifold $P$ of $G$. We also
give several equivalent representations of the submanifold $P$.
Some properties of the set $gK\cap P$ are also discussed, where
$gK$
is a coset space in $G$.\\


\footnotetext[1]{This work is supported by the 973 Project
Foundation of China ($\sharp$TG1999075102).}


\vskip 1.0cm
\section{Introduction}

\vskip 0.3cm Suppose $G$ is a connected Lie group, $K$ a closed
subgroup of $G$. Then $G/K$ is a homogeneous space. When $G$,
being considered as a principal $K$-bundle, is trivial, there is a
global section $S$ of this bundle. In this case, there is a
natural isomorphism $G/K\cong S$. This happens when $G$ is
semisimple and $K$ is a maximal compact subgroup of $G$, thanks to
the Cartan decomposition. If $G$ is not a trivial $K$-bundle, we
seldom consider $G/K$ as a submanifold of $G$. In this paper, we
show that when $G/K$ has the structure of Riemannian symmetric
space where $K$ is the closed subgroup of $G$ consist of the fixed
points of the involution of $G$ defining the symmetric space, we
can embed $G/K$ in $G$ as a closed submanifold in a good manner.

\vskip 0.3cm More precisely, let $G$ be a connected Lie group,
$\sigma$ an involution of $G$. Let $K=\{g\in G|\sigma(g)=g\}$,
then $K$ is a closed subgroup of $G$. We suppose there exists a
$G$-invariant Riemannian structure on $G/K$. Then $G/K$ becomes a
Riemannian (globally) symmetric space (see Helgason \cite{He}). In
this case, we call the triple $(G,\sigma, K)$ a {\it Riemannian
symmetric triple}. The differential $d\sigma$ of $\sigma$ gives an
involution of $\g$, the Lie algebra of $G$. Let $\ak$ and $\ap$ be
the eigenspaces of $d\sigma$ in $\g$ with eigenvalues $1$ and
$-1$, respectively. Then $\g=\ak\oplus\ap$. They satisfy the
relations $ [\ak,\ak]\subset\ak$, $[\ak,\ap]\subset\ap$, and
$[\ap,\ap]\subset\ak$. Note that $\ak$ is the Lie algebra of $K$.
Note also that there exists a $G$-invariant Riemannian structure
on $G/K$ if and only if $\Ad(K)|_\ap=\{\Ad(k)|_\ap\in GL(\ap)|k\in
K\}$ is compact.

\vskip 0.3cm In \S 2, we will prove that $P=\exp(\ap)$ is a closed
submanifold of $G$, and there is a natural isomorphism $G/K\cong
P$. That is, $G/K$ can be embedded as a closed submanifold of $G$.
The most hard part of the proof is the closedness of $P$. We will
deduce it by proving the fact that $P$ coincides with the
connected component $R_0$ containing $e$ of the set $R=\{g\in
G|\sigma(g)=g^{-1}\}$. In fact, we will give several equivalent
representations of $P$, namely $P=Q=R_0=R'_0=R^2$. This is our
main Theorem 2.5 in \S 2. Then, as a corollary, we get the natural
embedding of $G/K$ in $G$. In the case that $G$ is semisimple, the
relation $P=R_0$ was mentioned in Hermann \cite{Her}, Ch. 6, but
the author did not give a proof there. It is interesting to notice
that we can prove each connected component of $R$ is a closed
submanifold of $G$, but different components may have different
dimensions.

\vskip 0.3cm Even if $P$ is a closed submanifold of $G$, it is not
a global section of the principal $K$-bundle $G\rightarrow G/K$ in
general. But since $\g=\ak\oplus\ap$, $P$ is a local section
around $[e]=K$. Then it is naturally to ask that how far is $P$
from being a global section. This will be discussed in \S3. We
will prove, among other things, that $gK\cap P\neq\emptyset$ for
each coset space $gK$, and almost all coset space $gK$ intersects
$P$ transversally.

\vskip 0.3cm We would like to thank Professor U. Magnea for
helpful discussion. In fact, his work with Professor M. Caselle on
random matrix theory of symmetric spaces \cite{CM} is the
motivation for us to concentrate on the problems of this paper.


\vskip 1.0cm
\section{Realization of symmetric spaces in Lie groups}

\vskip 0.3cm We always suppose $(G,\sigma, K)$ is a Riemannian
symmetric triple as we have defined in section \S1 from now on. We
construct the sets
\[P=\exp(\ap),\]
\[Q=\{g\sigma(g)^{-1}|g\in G\},\]
\[R=\{g\in G|\sigma(g)=g^{-1}\}.\]
Let $R_0$, $R'_0$ be the connected component and the path
component of $R$ containing the identity $e$, respectively, and
let $R^2=\{g^2|g\in R\}$. Let $K_0$ be the identity component of
$K$.

\begin{lemma}\label{2.1}
The map $\Phi:P\times K_0\rightarrow G, \Phi(p,k)=pk$ is
surjective.
\end{lemma}

\begin{proof}
$\forall g\in G$, we prove that there are $k\in K_0$ and $X\in\ap$
such that $g=e^Xk$. Note that $G/K_0$ is also a symmetric space,
which is a covering space of $G/K$. We denote $x_0=[e]\in G/K_0$,
and let $x_1=gx_0$. Let $\gamma(t)$ be a geodesic in $G/K_0$ such
that $\gamma(0)=x_0$ and $\gamma(1)=x_1$. Then $\gamma(t)$ is of
the form $\gamma(t)=e^{tX}x_0$ for some $X\in\ap$ (see Helgason
\cite{He}, Ch. IV, \S3). So $x_1=e^Xx_0$. Let $k=e^{-X}g$, we have
$kx_0=e^{-X}(gx_0)=e^{-X}x_1=x_0$. So $k\in K_0$, and $g=e^Xk$.
\end{proof}

\begin{lemma}\label{2.2}
For all $p,p'\in P$, we have $pp'p\in P$.
\end{lemma}

\begin{proof}
Suppose $p'=e^X$, where $X\in\ap$. Let $p_1=e^{X/2}\in P$, then
$p'=p_1^2$. By Lemma \ref{2.1}, $pp_1=p_2k$ for some $p_2\in P,
k\in K_0$. Then we have
$p^{-1}p_1^{-1}=\sigma(pp_1)=\sigma(p_2k)=p_2^{-1}k$, this implies
$p_1p=k^{-1}p_2$. So $pp'p=(pp_1)(p_1p)=(p_2k)(k^{-1}p_2)=p_2^2\in
P$.
\end{proof}

A neighborhood $U$ of $0$ in $\g$ is {\it symmetric} if $X\in U$
implies $-X\in U$.

\begin{lemma}\label{2.4}
Suppose $U$ is a symmetric neighborhood of $0$ in $\g$ with
$d\sigma(U)=U$ such that $\exp|_U:U\rightarrow \exp(U)$ is a
diffeomorphism. If $g\in\exp(U)$ satisfies $\sigma(g)=g^{-1}$,
then $g\in P$.
\end{lemma}

\begin{proof}
Suppose $g=e^X$, $X\in U$. Applying the equation $\exp\circ
d\sigma=\sigma\circ\exp$ to $X$, we have
$\exp(d\sigma(X))=\sigma(\exp(X))=\exp(-X)$. Since $U$ is
symmetric and $d\sigma(U)=U$, $d\sigma(X),-X\in U$. But $\exp$ is
injective on $U$, so we have $d\sigma(X)=-X$. This implies
$X\in\ap$, so $g=e^X\in P$.
\end{proof}

\begin{lemma}\label{2.5}
Suppose $U$ is a symmetric neighborhood of $0$ in $\g$ with
$d\sigma(U)=U$ such that $\exp|_U$ is a diffeomorphism onto its
image, and suppose $p\in P$. If $g\in p\exp(U)p$ satisfies
$\sigma(g)=g^{-1}$. Then $g\in P$.
\end{lemma}

\begin{proof}
Since $g\in p\exp(U)p$, $p^{-1}gp^{-1}\in\exp(U)$. But
$\sigma(p^{-1}gp^{-1})=pg^{-1}p=(p^{-1}gp^{-1})^{-1}$. By Lemma
\ref{2.4}, $p^{-1}gp^{-1}\in P$. Then by Lemma \ref{2.2}, $g\in
P$.
\end{proof}

Now we are prepared to formulate our main Theorem in this section.

\begin{theorem}\label{a}
Suppose $(G,\sigma, K)$ is a Riemannian symmetric triple, and let
the subsets $P,Q,R,R_0,R'_0,R^2$ of $G$ be as defined above. Then
$P=Q=R_0=R'_0=R^2$.
\end{theorem}

\begin{proof}
We prove $P\subset R^2\subset Q\subset R'_0\subset P$ and
$R_0=R'_0$.

\vskip 0.3cm (i) ``$P\subset R^2.$" Suppose $g\in P$, then $g=e^X$
for some $X\in\ap$. But
$\sigma(e^{X/2})=e^{d\sigma(X/2)}=e^{-X/2}$, so $e^{X/2}\in R$,
and then $g=(e^{X/2})^2\in R^2$.

\vskip 0.3cm (ii) ``$R^2\subset Q.$" Suppose $g\in R^2$, then
$g=h^2$, $\sigma(h)=h^{-1}$. Now $h\sigma(h)^{-1}=h^2=g$, so $g\in
Q$.

\vskip 0.3cm (iii) ``$Q\subset R'_0.$" For $g\sigma(g)^{-1}\in Q$,
$\sigma(g\sigma(g)^{-1})=\sigma(g)g^{-1}=(g\sigma(g)^{-1})^{-1}$,
so $Q\subset R$. But $Q$ is path connected and containing $e$, so
$Q\subset R'_0$.

\vskip 0.3cm (iv) ``$R'_0\subset P.$" We first suppose that $G$ is
simply connected. Let $g:[0,1]\rightarrow R'_0$ be a continuous
path in $R'_0$ with $g(0)=e$, it suffices to prove $g(1)\in P$.
Let $S=\{t\in[0,1]|g(t)\in P\}$. Since $g(0)\in P$, $0\in S$, so
$S\neq\emptyset$. We will prove that $S$ is open and closed. Then
by the connectedness of $[0,1]$, $S=[0,1]$, and then we will have
$g(1)\in P$.

\vskip 0.3cm For the openness of $S$, suppose $t_0\in S$, that is
$g(t_0)=e^X$ for some $X\in\ap$. Let $p=e^{X/2}$, then
$g(t_0)=p^2$. Let $U$ be a symmetric neighborhood of $0$ in $\g$
with $d\sigma(U)=U$ such that $\exp|_U$ is a diffeomorphism onto
its image. Then $p\exp(U)p$ is a neighborhood of $g(t_0)$. So
there is an open neighborhood $(t_1,t_2)$ of $t_0$ such that
$g(t)\in p\exp(U)p, \forall t\in(t_1,t_2)$. By Lemma \ref{2.5},
$g(t)\in P$, that is $t\in S, \forall t\in(t_1,t_2)$. This proves
the openness.

\vskip 0.3cm To prove that $S$ is closed, we endow a left
invariant Riemannian structure on $G$. This induces a left
invariant metric $d(\cdot,\cdot)$ on $G$. Since $G$ is simply
connected, there is an $\Ad(G)$ invariant, $d\sigma$-invariant
symmetric neighborhood $V$ of $0$ in $\g$ such that $\exp|_V$ is a
diffeomorphism onto its image (see Varadarajan \cite{Va}, Theorem
2.14.6). Then $g\exp(V)g^{-1}=\exp(V), \forall g\in G$. Let $r>0$
such that $B_r(e)=\{g\in G|d(e,g)<r\}\subset \exp(V)$. Suppose
$\{t_n\}_{n\in\mathbb{N}}\subset S$ is a sequence such that
$\lim_{n\rightarrow\infty}t_n=t_0$, we prove $t_0\in S$. Choose
$N\in\mathbb{N}$ such that $d(g(t_N),g(t_0))<r$. Since $g(t_N)\in
P$, $g(t_N)=p^2$ for some $p\in P$. So $g(t_0)\in
g(t_N)B_r(e)=p^2B_r(e)\subset
p^2\exp(V)=p^2(p^{-1}\exp(V)p)=p\exp(V)p$. By Lemma \ref{2.5},
$g(t_0)\in P$. Hence $S$ is closed. This conclude the proof when
$G$ is simply connected.

\vskip 0.3cm For general $G$, let $\widetilde{G}$ be its universal
covering group with covering map $\pi:\widetilde{G}\rightarrow G$.
Let the corresponding involution of $\widetilde{G}$ is
$\widetilde{\sigma}$, and let $\widetilde{R}'_0$, $\widetilde{P}$
be the corresponding subsets of $\widetilde{G}$. We claim that
$R'_0\subset\pi(\widetilde{R}'_0)$. In fact, suppose $g\in R'_0$.
Then there is a continuous path $g(t) (t\in[0,1])$ in $R'_0$ such
that $g(0)=e$ and $g(1)=g$. Let $\widetilde{g}(t)$ be a lift of
$g(t)$ to $\widetilde{G}$ such that $\widetilde{g}(0)=e$. Then
$\pi(\widetilde{g}(t)\widetilde{\sigma}(\widetilde{g}(t)))=g(t)\sigma(g(t))=e$,
that is
$\widetilde{g}(t)\widetilde{\sigma}(\widetilde{g}(t))\in\ker(\pi)$.
But $\ker(\pi)$ is discrete and
$\widetilde{g}(0)\widetilde{\sigma}(\widetilde{g}(0))=e$, so
$\widetilde{g}(t)\widetilde{\sigma}(\widetilde{g}(t))=e$, that is
$\widetilde{g}(t)\in \widetilde{R}'_0$. In particular,
$\widetilde{g}(1)\in \widetilde{R}'_0$. But
$g=g(1)=\pi(\widetilde{g}(1))$, so $g\in\pi(\widetilde{R}'_0)$.
Hence we have
$R'_0\subset\pi(\widetilde{R}'_0)\subset\pi(\widetilde{P})=P$.
Then (iv) is proved.

\vskip 0.3cm (v) ``$R_0=R'_0.$" It is well known that $R'_0\subset
R_0$. To prove the converse, we let
$V=\{X\in\g||\mathrm{Im}(\lambda)|<\pi \text{ for each eigenvalue
} \lambda \text{ of } \ad(X)\}$. Then $V$ is an $\Ad(G)$
invariant, $d\sigma$-invariant symmetric neighborhood of $0\in\g$.
By Varadarajan \cite{Va}, Theorem 2.14.6, there is a discrete
additive subgroup $\Gamma$ of $\g$ such that for $X, X'\in V$,
$e^X=e^{X'}$ if and only if $X-X'\in\Gamma$. Choose a neighborhood
$U\subset V$ of $0\in\g$ and $r_0>0$ such that $\exp|_U$ is a
diffeomorphism onto its image and $\exp(U)=B_{r_0}(e)$. ($G$ being
endowed a left invariant Riemannian structure.) For a continuous
function $\rho:G\rightarrow(0,r_0)$, let $N_\rho=\bigcup_{g\in
R'_0}B_{\rho(g)}(g)$. Then $N_\rho$ is an open neighborhood of
$R'_0$. It is easy to prove that $R'_0\subset
N_{\rho/2}\subset\overline{N_{\rho/2}}\subset N_\rho$. We will
prove that for sufficient small $\rho:G\rightarrow(0,r_0)$,
$N_\rho\cap R=R'_0$. So $R'_0$ is open in $R$. $N_\rho\cap R=R'_0$
implies $\overline{N_{\rho/2}}\cap R=R'_0$, so $R'_0$ is closed in
$R$. This means that $R'_0$ is in fact a connected component of
$R$. So we will have $R_0=R'_0$.

\vskip 0.3cm Now we prove $N_\rho\cap R=R'_0$ for sufficient small
$\rho$. Since $\Gamma\subset\g$ is discrete, there is $\epsilon>0$
such that $B_\epsilon(0)\cap\Gamma=\{0\}$. Let $K_1\subset
K_2\subset\cdots\subset K_n\subset\cdots$ be a sequence of compact
subsets of $G$ such that $\bigcup_{n=1}^\infty K_n=G$. Let
$C_n=\sup_{g\in K_n}\|d\sigma+\Ad(g)\|+1$. Choose $r_n\in(0,r_0)$
such that $B_{r_n}(e)\subset\exp(B_{\epsilon/C_n}(0))$. We claim
that if the function $\rho:G\rightarrow(0,r_0)$ satisfies
$\rho(g)<r_n$, for $g\in K_n, \forall n\in\mathbb{N}$, then
$N_\rho\cap R=R'_0$. In fact, for $g'\in N_\rho\cap R$, by the
definition of $N_\rho$, there exists $g\in R'_0$ such that $g'\in
B_{\rho(g)}(g)$. Suppose $g\in K_n$. Let $g'=gh$. By the left
invariance of the Riemannian structure, $h\in
B_{\rho(g)}(e)\subset B_{r_n}(e)$. Suppose $h=e^X$, where $X\in
U$. By $B_{r_n}(e)\subset\exp(B_{\epsilon/C_n}(0))$, we know that
$|X|<\epsilon/C_n$. Since $g,g'\in R$,
$g^{-1}\sigma(h)=\sigma(gh)=\sigma(g')=g'^{-1}=h^{-1}g^{-1}$, that
is $\sigma(h)=gh^{-1}g^{-1}$. But $h=e^X$, this implies
$\exp(d\sigma(X))=\exp(-\Ad(g)X)$. Since $d\sigma(X),-\Ad(g)X\in
V$, $d\sigma(X)-(-\Ad(g)X)=(d\sigma+\Ad(g))X\in\Gamma$. But
$|(d\sigma+\Ad(g))X|\leq\|d\sigma+\Ad(g)\|\cdot|X|<C_n\cdot\epsilon/C_n=\varepsilon$.
So $(d\sigma+\Ad(g))X\in B_\epsilon(0)\cap\Gamma=\{0\}$, that is,
$d\sigma(X)=-\Ad(g)X$. Now let $\gamma(t)=ge^{tX}$, then
$\gamma(0)=g$, $\gamma(1)=g'$. But
$\sigma(\gamma(t))=\sigma(ge^{tX})=g^{-1}\exp(td\sigma(X))=g^{-1}\exp(-t\Ad(g)X)
=g^{-1}ge^{-tX}g^{-1}=(ge^{tX})^{-1}=\gamma(t)^{-1}.$ So
$\gamma(t)\in R$. This proves $g'=\gamma(1)\in R'_0$. (v) is
proved.
\end{proof}

\begin{remark}
In general, $P\subsetneqq R$, even in some very simple cases. For
example, Let $G=SL(2n,\RR)$, $\sigma(g)=(g^t)^{-1}$. Then
$\diag(-1,\cdots,-1)$ is in $R$, but not in $P$. It is obvious by
the above theorem that $P=R$ if and only if $R$ is connected (or
path connected).
\end{remark}

\begin{remark}
Using the same method as in the proof of (v) of the above theorem,
we can prove that for each path component $R'_i$ of $R$, there are
open subsets $U_i, V_i$ of $G$ such that $R'_i\subset
V_i\subset\overline{V_i}\subset U_i$ and $U_i\cap R=R'_i$. Thus
$R'_i$ is open and closed in $R$, and then $R'_i$ is in fact a
connected component of $R$. The proof of (v) of the above theorem
also shows that $R'_0$ is a closed submanifold of $G$ of dimension
$\dim(\ker(d\sigma+\Ad(g))),\; g\in R'_0$. Similarly, each $R'_i$
is a closed submanifold of $G$ of dimension
$\dim(\ker(d\sigma+\Ad(g))),\; g\in R'_i$. So each connected
component of $R$ is a closed submanifold of $G$. But different
components of $R$ can have different dimensions. For example, let
$G=SO(5)$ and $\sigma(g)=sgs$, where $s=\diag(1,-1,-1,-1,-1)$.
Then $\dim R'_0=\dim(\ker(d\sigma+\Ad(e)))=4$. But
$g_0=\diag(-1,-1,-1,-1,1)\in R$, and the connected component of
$R$ containing $g_0$ has dimension
$\dim(\ker(d\sigma+\Ad(g_0)))=6$. We leave the detail of the
proofs of these conclusions to the reader.
\end{remark}

We define the {\it twisted conjugate action} of $G$ on $G$ by
$\tau_g(h)=gh\sigma(g)^{-1}$. Then $Q=P$ is the orbit of this
action containing the identity $e$. The next conclusion says that
the symmetric space $G/K$ can be embedded in $G$ as a closed
submanifold, which is just the set $P\subset G$.

\begin{corollary}\label{b}
$P$ is a closed submanifold of $G$. The map
$\varphi:G/K\rightarrow P$ defined by
$\varphi([g])=g\sigma(g)^{-1}$ is a diffeomorphism. Under the
actions of $G$ by left multiplication on $G/K$ and by the twisted
conjugate action on $P$, the isomorphism $\varphi$ is equivariant.
\end{corollary}

\begin{proof}
The twisted conjugate action $\tau$ is smooth, so $Q$, as an orbit
of this action, is an immersion submanifold of $G$. As a connected
component of $R$, $R_0$ is closed in $R$. But $R$ is closed in
$G$, so $R_0$ is closed in $G$. By Theorem \ref{a}, $P=Q=R_0$ is a
closed submanifold of $G$. Notice that the isotropic subgroup of
the action $\tau$ associated with the identity $e$ is just $K$, so
$\varphi:G/K\rightarrow P$ is a diffeomorphism. It is obviously
equivariant.
\end{proof}

\begin{remark}
Corollary \ref{b} says that the symmetric space $G/K$ can be
realized in $G$ as a closed submanifold, which is just $P$. But we
should point out that for any subgroup $K'$ of $G$ satisfying
$K_0\subset K'\subset K$, $G/K'$ is also a symmetric space. In
general, $G/K'$ can not be embedded in $G$ as closed submanifold.
The submanifold $P\subset G$ is isomorphic at most one of such
$G/K'$, and this happens if and only if $K'=K$.
\end{remark}

\begin{remark}
If $G$ is compact, the closedness of $P$ can be implied from the
fact that $P=R^2$, which is much easier to be obtained than
$P=R_0$. In fact, to prove that $P$ is a closed submanifold when
$G$ is compact, we need only to show that $P=R^2=Q$. $P\subset
R^2\subset Q$ have been showed in the proof of Theorem \ref{a},
(i) and (ii). $Q\subset P$ can be proved by the following simple
argument. Let $g=h\sigma(h)^{-1}\in Q$. Since $h$ can be expressed
as $h=pk$, where $p\in P, k\in K$,
$g=(pk)(\sigma(pk))^{-1}=(pk)(p^{-1}k)^{-1}=p^2\in P$.
\end{remark}


\vskip 1.0cm
\section{How far is $P$ from being a section}

\vskip 0.3cm We regard $G$ as a principal $K$-bundle with base
space $G/K$. We have proved that $P=\exp(\ap)$ is a closed
submanifold of $G$. It is obvious that the tangent space $T_eP$ of
$P$ at $e$ is $\ap$. But $\g=\ak\oplus\ap$. So $P$ is a local
section of the bundle $G\rightarrow G/K$ around $[e]$. If $G$ is
semisimple (or more generally, connected reductive, in the sense
of Knapp \cite{Kn}) and $\sigma$ is the global Cartan involution
of $G$, then $K=K_0$, and by the Cartan decomposition, the map
$\Phi:P\times K\rightarrow G$ defined in Lemma \ref{2.1} is a
diffeomorphism. So $G$ is a trivial $K$-bundle, and $P$ is a
global section. This means that for each coset space $gK$, $gK\cap
P$ consists of just one point. But in general $G$ is not a trivial
$K$-bundle. So we may ask the question: for a coset space $gK$,
how many points in $gK\cap P$?

\begin{theorem}
For each coset space $gK_0$, there is a homeomorphism between
$gK_0\cap P$ and $\Phi^{-1}(g)$. In particular, $gK_0\cap
P\neq\emptyset$.
\end{theorem}

\begin{proof}
Let $\pi_1:P\times K_0\rightarrow P$ be the projection to the
first factor. We prove that $\pi_1|_{\Phi^{-1}(g)}$ is a
homeomorphism between $\Phi^{-1}(g)$ and $gK_0\cap P$. First, let
$(p,k)\in\Phi^{-1}(g)$, then $g=pk$. So $\pi_1(p,k)=p=gk^{-1}\in
gK_0\cap P$. This proves $\pi_1(\Phi^{-1}(g))\subset gK_0\cap P$.
Let $(p_1,k_1), (p_2,k_2)\in\Phi^{-1}(g)$ and
$(p_1,k_1)\neq(p_2,k_2)$. Since $p_1k_1=p_2k_2=g$, $p_1\neq p_2$.
This shows $\pi_1|_{\Phi^{-1}(g)}$ is injective. Let $p\in
gK_0\cap P$, then there is some $k\in K_0$ such that $p=gk$. So
$(p,k^{-1})\in\Phi^{-1}(g)$, and $\pi_1(p,k^{-1})=p$. This means
that $\pi_1|_{\Phi^{-1}(g)}:\Phi^{-1}(g)\rightarrow gK_0\cap P$ is
surjective. Since $\pi_1$ is continuous and open, so is
$\pi_1|_{\Phi^{-1}(g)}$. Hence
$\pi_1|_{\Phi^{-1}(g)}:\Phi^{-1}(g)\rightarrow gK_0\cap P$ is a
homeomorphism. By Lemma \ref{2.1}, $\Phi$ is surjective. So
$gK_0\cap P\cong\Phi^{-1}(g)\neq\emptyset$. This proves the
theorem.
\end{proof}

\begin{corollary}
For each coset space $gK$, $gK\cap P\neq\emptyset$.\qed
\end{corollary}

Similar to the map $\Phi:P\times K_0\rightarrow G$, we can define
the map $\Phi':P\times K\rightarrow G$ by $\Phi(p,k)=pk$. It is
easy to see that $\Phi'$ satisfies all the properties of $\Phi$
that we have mentioned above.

\vskip 0.3cm In the following we denote the left and right
translations of $g\in G$ by $L_g$ and $R_g$, respectively.

\begin{lemma}\label{3.3}
Let $g\in G$. Then $g$ is a regular value of $\Phi'$ if and only
if $gK$ intersects $P$ transversally.
\end{lemma}

\begin{proof}
Suppose $(p,k)\in\Phi'^{-1}(g)$. Then
\begin{align*}
&\im (d\Phi')_{(p,k)}\\
=&(d\Phi')_{(p,k)}(T_{(p,k)}(P\times\{k\}))+(d\Phi')_{(p,k)}(T_{(p,k)}(\{p\}\times
K))\\
=&T_g(Pk)+T_g(pK)\\
=&(dR_k)_p(T_pP)+(dR_k)_p(T_p(gK))\\
=&(dR_k)_p(T_pP+T_p(gK)).
\end{align*}
Since $(dR_k)_p$ is an isomorphism and $\pi_1(\Phi'^{-1}(g))=
gK\cap P$, we have
\begin{align*}
&g \;\text {is a regular value of}\; \Phi'\\
\Leftrightarrow &\im (d\Phi')_{(p,k)}=T_gG,\quad\forall(p,k)\in\Phi'^{-1}(g)\\
\Leftrightarrow &T_pP+T_p(gK)=T_pG,\quad\forall(p,k)\in\Phi'^{-1}(g)\\
\Leftrightarrow &T_pP+T_p(gK)=T_pG,\quad\forall p\in gK\cap P\\
\Leftrightarrow &gK \;\text {intersects}\; P \;\text
{transversally}.
\end{align*}
\end{proof}

\begin{lemma}\label{3.4}
The set of all regular values of $\Phi'$ is right $K$-invariant.
\end{lemma}

\begin{proof}
For $g\in G$ and $k_1\in K$,
$(p,k)\in\Phi'^{-1}(g)\Leftrightarrow(p,kk_1)\in\Phi'^{-1}(gk_1)$.
Since $$R_{k_1}\circ\Phi'=\Phi'\circ(id\times R_{k_1}),$$ for
$(p,k)\in\Phi'^{-1}(g)$, we have
$$(dR_{k_1})_g\circ(d\Phi')_{(p,k)}=(d\Phi')_{(p,kk_1)}\circ(id\times dR_{k_1})_{(p,k)}.$$
Since $(dR_{k_1})_g$ and $(id\times dR_{k_1})_{(p,k)}$ are
isomorphisms, $(d\Phi')_{(p,k)}$ is an isomorphism if and only if
$(d\Phi')_{(p,kk_1)}$ is an isomorphism. So $g$ is a regular value
if and only if $gk_1$ is a regular value. This proves the lemma.
\end{proof}

\begin{theorem}\label{3.5}
For almost all coset space $gK$ in $G/K$, $gK$ intersects $P$
transversally.
\end{theorem}

\begin{proof}
Let $G_r$ be the set of all regular values of $\Phi'$. By Sard's
Theorem, $G\setminus G_r$ is a set with measure zero. But Lemma
\ref{3.4} tells us that $G_r$ is the union of some coset spaces
$gK$. So by choosing local trivializations of the principal bundle
$\pi:G\rightarrow G/K$ and using Fubini's Theorem, we know that
$\pi(G\setminus G_r)=(G/K)\setminus\pi(G_r)$ has measure zero in
$G/K$. By Lemma \ref{3.3}, $[g]=gK$ intersects $P$ transversally,
$\forall\; [g]\in\pi(G_r)$. This proves the theorem.
\end{proof}

\begin{corollary}
For almost all coset space $gK$ in $G/K$, $gK\cap P$ is a discrete
set. In particular, if $K$ is compact, then $gK\cap P$ is a finite
set for almost all coset space $gK$ in $G/K$.
\end{corollary}

\begin{proof}
Since $\dim(gK) +\dim P=\dim G$, $gK$ intersects $P$ transversally
implies that $gK\cap P$ is a 0-dimensional submanifold of $G$,
which is discrete. In particular, if $K$ is compact, then so is
$gK$. But $gK\cap P\subset gK$. So $gK$ intersects $P$
transversally implies $gK\cap P$ is a finite set. By Theorem
\ref{3.5}, the corollary holds.
\end{proof}

Denote $x_0=[e]\in G/K$. We know that a curve $\gamma(t)
(t\in\RR)$ in $G/K$ with $\gamma(0)=x_0$ is a geodesic if and only
if $\gamma(t)=e^{tX}x_0$ for some $X\in\ap$, and such $X$ is
unique.

\begin{theorem}\label{3.7}
Let $g\in G$. If for every two geodesics $\gamma_i(t)=e^{tX_i}x_0
(i=1,2)$ through the point $[g]$, where $X_i\in\ap$, we have
$[X_1,X_2]=0$. Then $\#(gK\cap P)\leq\#(K\cap P)$.
\end{theorem}

\begin{proof}
Suppose $g=pk$, where $p\in P, k\in K$. We prove
$L_{p^{-1}}(gK\cap P)\subset(K\cap P)$. Suppose $p'\in(gK\cap P)$.
Since $p\in gK$, $L_{p^{-1}}(p')=p^{-1}p'\in K$. Let $p=e^X,
p'=e^{X'}$, where $X,X'\in\ap$. Since $[p]=[p']=[g]$, the two
geodesics $\gamma_1(t)=e^{tX}x_0$ and $\gamma_2(t)=e^{tX'}x_0$
satisfy $\gamma_1(1)=\gamma_2(1)=[p]$. By the conditions of the
theorem, $[X,X']=0$. So $L_{p^{-1}}(p')=e^{-X}e^{X'}=e^{X'-X}\in
P$. This proves $L_{p^{-1}}(gK\cap P)\subset(K\cap P)$. But
$L_{p^{-1}}$ is injective, this proves the theorem.
\end{proof}

By Lemma \ref{2.1}, all coset space $gK$ is of the form $pK$ for
some $p\in P$. We conclude this section by an example in which we
show what the set $pK\cap P$ is for each $p\in P$.

\begin{example}
Let $G=SU(2)$, $\sigma(g)=(g^t)^{-1}$. Then $K=SO(2)$.
$$R=\{g\in SU(2)|g^t=g\}
=\Big\{
\begin{pmatrix}
a+bi & ci\\
ci & a-bi
\end{pmatrix}
\Big|a,b,c\in\RR, a^2+b^2+c^2=1\Big\}.$$ So $R\cong S^2$ is
connected, and then $P=R$. For $p\in P$, we show what the set
$pK\cap P$ is. The element of $pK$ has the form $pk, k\in K$. But
$$pk\in P\Leftrightarrow \sigma(pk)=(pk)^{-1}\Leftrightarrow p^{-1}k=k^{-1}p^{-1}\Leftrightarrow kp=pk^{-1}.$$
Let $k=\begin{pmatrix}
\cos \theta & -\sin \theta\\
\sin \theta & \cos \theta
\end{pmatrix}$, $p=\begin{pmatrix}
a+bi & ci\\
ci & a-bi
\end{pmatrix}$. Then it is easy to show that $kp=pk^{-1}\Leftrightarrow a\sin
\theta=0$. So if $a\neq0$, $pk\in P\Leftrightarrow \sin
\theta=0\Leftrightarrow k=\pm I$. In this case $pK\cap P=\{\pm
p\}$. In particular, $K\cap P=\{\pm I\}$. If $a=0$, then $\forall
k\in K, pk\in P$. This implies $pK\subset P$. So in this case,
$pK\cap P=pK$. It should be noted that all $p=\begin{pmatrix}
bi & ci\\
ci & -bi
\end{pmatrix}$ correspond the same coset space $pK$, which is the
antipodal point of $[e]$ in the symmetric space $SU(2)/SO(2)\cong
S^2$.\qed
\end{example}

\end{document}